\documentstyle[12pt,amscd]{amsart}
\newtheorem{Theorem}{Theorem}[section]
\newtheorem{Lemma}[Theorem]{Lemma}
\newtheorem{Corollary}[Theorem]{Corollary}

\newtheorem{Remark}[Theorem]{Remark}
\newtheorem{Definition}[Theorem]{Definition}

\def\reg{\operatorname{reg}}

\def\sk{\smallskip\par}
\def\mm{{\mathfrak m}}
\def\aa{{\mathfrak a}}

\def\qq{{\mathfrak q}}

\begin{document}
\title{Castelnuovo-Mumford regularity of \\
 associated graded modules  \\ of  $k$-Buchsbaum modules}
\maketitle
 \begin{center}
 LE XUAN DUNG\\
 Faculty of  Natural Sciences,  Hong Duc University\\
 No 565 Quang Trung, Dong Ve, Thanh Hoa, Vietnam\\
 E-mail: lexuandung@@hdu.edu.vn\\ [15pt]
\end{center}

 \begin{abstract}  
Bounds on the Castelnuovo-Mumford regularity of the
 associated graded modules  of $k$-Buchsbaum modules $M$ are given in terms of  $k$ and some other invariants of $M$.
\end{abstract}

\date{}
\section*{Introduction} \sk

Let $(A,\mm)$ be a commutative Noetherian local ring, $\aa$ an $\mm$-primary ideal and $M$ a finitely generated $A$-module.  The associated graded module to the $M$ with respect to $\aa$ is defined by $G_\aa(M) = \oplus_{n\ge 0}\aa^nM/\aa^{n+1}M$. The  Castelnuovo-Mumford regularity $\reg(G_\aa(M))$ of $G_\aa(M)$ can be used to study other invariants of $M$. Therefore bounding this invariant  is an important problem.  A breakthrough solution  was  given by  Rossi-Trung-Valla  for the case $M=A$ and $\aa= \mm$  in \cite{RTV}, which was then extended to other cases in \cite{L} and  \cite{DH1}, where the bound  is given in terms of  the so-called extended degree $D(\aa,M)$. Other bounds on Castelnuovo-Mumford regularity $\reg(G_\aa(M))$  are given in  \cite{DH2,DH3} in terms of Hilbert coefficients. Note that the so-called homological degree is an extended degree, see \cite[Definition 9.4.1]{Va} for the graded case and \cite{RTV, L} for the local case. In the case of  generalized Cohen-Macaulay modules, the homological degree can be expressed in terms of the so-called  Buchsbaum  invariant $I(M)$. 

In this paper, we restrict to the case of  generalized Cohen-Macaulay modules. Then a generalized Cohen-Macaulay module is a $k$-Buchsbaum module for some $k\ge 0$. In order to compute $I(M)$ one has to know all local cohomology modules $H^i_{\mm}(M)$ of $M$, where $i<\dim(M)$. This is of course not an easy problem.  The aim of this paper is to give a bound on $\reg(G_\aa(M))$ in terms of $k$ and two other invariants: the so-called reduction number and the length of a quotient module, which are much easier to compute. 

The paper is divided into two sections. In Section 1, we start with a few preliminary results on $k$-Buchsbaum modules. In Section 2, we give new bounds on $\reg(G_\aa(M))$   in terms of  $k$ (see Theorem \ref{a3.1.2}, Theorem \ref{a3.1.22} and Theorem \ref{main}).

\section{Preliminaries} 

Let $A$ be  a local ring with the maximal ideal $\mm$ and  $M$ be a finitely generated  $A$-module such that  $d = \dim A = \dim M \ge 1$.  We always assume that the residue field $A/\mm$ is infinite.

\begin{Definition} \label{a1} {\rm (See, e. g., \cite{HV})}
Let $k \ge 0$ be an integer and $1 \le t \le d$. $M$ is called a $(k,t)$-Buchsbaum module if for every s.o.p. $x_1,x_2,...,x_d$ of $M$ we have $$ \mm^kH_{\mm}^i(M/(x_1,...,x_j)M) = 0, $$
for all non-negative integers $i,j$  with $j \le t-1$ and $i+j < d$.

$A$ is called a $(k,t)$-Buchsbaum  ring if it is a $(k,t)$-Buchsbaum module over itself. A $(k,1)$-Buchsbaum module (ring) is also called a $k$-Buchsbaum module (resp. ring).
\end{Definition}
\begin{Remark}  \label{a2} 
$(0,t)$-, $(1,d)$-, and $1$-Buchsbaum module are exactly the Cohen-Macaulay, Buchsbaum, and quasi-Buchsbaum modules, resp. 
\end{Remark}

A module $M$ is called a {\it generalized Cohen-Macaulay module} if   $\ell(H_{\mm}^i(M)) < \infty$ for $i=0,1,...,d-1$, where $\ell$ denotes the length and $H_{\mm}^i(M) $ the $i$th local cohomology module of $M$ with respect to $\mm$ (see \cite{CST,T1}).

For a parameter ideal $Q$ of $M$, let $$I(Q,M):= \ell(M/ QM) - e(Q,M),$$
where $e(Q,M)$ denotes the multiplicity of $A$ with respect to $Q$,  and let
$$ I(M): = \sup \{ I(Q,M) | Q \text{ is a  parameter ideal of } M \}. $$
This number is finite thank to the following fundamental result in the theory of  generalized Cohen-Macaulay modules:
\begin{Lemma} {\rm \cite[(3.7)]{CST}}\label{a1.1.4}
Let  $M$ be a generalized Cohen-Macaulay module. Then
$$I(M) = \sum_{i=0}^{d-1}\binom{d-1}{i}\ell(H_{\mm}^i(M)).$$
\end{Lemma}
The notion of generalized Cohen-Macaulay modules extends that of Buchsbaum modules. Therefore the number $I(M)$ is some times called the {\it Buchsbaum invariant} of $M$.

Generalized Cohen-Macaulay modules may be characterized in different ways.  Recall that a s.o.p. $x_1,...,x_d$  of $M$ is called {\it standard} if 
$$ (x_1,...,x_d)H_{\mm}^{i}(M/(x_1,...,x_j)M)=0$$
for all non-negative integer $i,j$ with $i+j < d$.

An ideal $I \subseteq \mm$ is called a {\it standard ideal} of $M$ if every s.o.p of $M$ contained in $I$ is a standard s.o.p of $M$. 

From (3.3) and (3.7) of \cite{CST}; Theorem 2.1, Theorem 2.5 and  Proposition 2.10 of \cite{T1} and \cite[Remark 1.3]{HV}, we have 

\begin{Lemma}  \label{GCM}  The following conditions  are equivalent:
\begin{itemize}
\item[{\rm (i)}]   $M$ is a generalized Cohen-Macaulay,
\item[{\rm (ii)}]  $I(M) < \infty $,
\item[{\rm (iii)}] There exists a positive integer $l$ such that $\mm^l$ is a standard ideal of $M$, 
 \item[{\rm (iv)}] $M$ is a $(k,t)$-Buchsbaum module for some $k \ge 0$ and for a certain $t$ with $1\le t\le d$.
\end{itemize}
\end{Lemma}
From the Definition \ref{a1} and \cite[Corollary 2.3]{T1}, we get
\begin{Lemma} \label{a3}
If $M$ is a $(k,d)$-Buchsbaum module, then $\mm^k$ is standard ideal of $M$ and $$ (x_1^{k},...,x_d^{k})M\cap H_{\mm}^0(M)=0,$$
for any s.o.p. $x_1,...,x_d$ of $M$.
\end{Lemma}

For a real number $\alpha $, let $\lceil \alpha \rceil$ denote the smallest integer $n$ such that $n\ge \alpha $.
\begin{Lemma}\label{a1.3.3}
Let $M$ be a $(k,d)$-Buchsbaum module with  $k>0$, $s$  a positive integer   and  $x_1,...,x_d \in \mm^s$ a s.o.p of $M.$  Let $k' = \lceil  \frac{k}{s} \rceil$. Then
$$(x_1,...,x_d)^{(k'-1)d+1}M\cap H_{\mm}^0(M)=0.$$
\end{Lemma}
\begin{pf}
Since $x_1,...,x_d \in \mm^s$, we  have $x_1^{k'},...,x_d^{k'} \in \mm^k $.
By Definition \ref{a1},  $\mm^{k}$ is a standard ideal of $M,$  hence  $ x_1^{k'},...,x_d^{k'} $ is a standard s.o.p of $M.$  By Lemma \ref{a3}, we get
$$ (x_1^{k'},...,x_d^{k'})M\cap H_{\mm}^0(M)=0.$$
Other hand $(x_1,...,x_d)^{(k'-1)d+1}\subseteq (x_1^{k'},...,x_d^{k'}).$
Therefore $$(x_1,...,x_d)^{(k'-1)d+1}M\cap H_{\mm}^0(M) =0 .$$
\end{pf}

Further we recall some more notions and notations. Let $R= \oplus_{n\ge 0}R_n$ be a Noetherian standard graded ring over a local Artinian ring $(R_0,\mm_0)$.  Let $E$ be a finitely generated graded $R$-module of dimension $d$.  For $0\le i \le d$, put
 $$ a_i(E) =
\sup \{n|\ H_{R_+}^i(E)_n \ne 0 \} ,$$
where $R_+ = \oplus _{n>  0} R_n$. The {\it Castelnuovo-Mumford regularity} of $E$  is defined by
$$\reg(E) = \max \{ a_i(E) + i \mid  0\le  i \leq  d  \},$$
and the {\it Castelnuovo-Mumford regularity of $E$  at and above level $1$},   is defined by
$$\reg^1(E) = \max \{ a_i(E) + i \mid  1\leq i \leq  d \}.$$

The {\it associated graded module} of $M$ 
 with respect to  an ideal $\aa$ of $A$ is defined by
  $$G_\aa(M) = \bigoplus_{n\geq 0}\aa^nM/\aa^{n+1}M. $$
This is a module over the associated graded ring $G =  \bigoplus_{n\geq 0}\aa^n/\aa^{n+1}.$  From now on, we assume that $\aa$ is an $\mm$-primary ideal. Then $A/\aa$ is an Artinian ring, so  $\reg(G_\aa(M)$ is well defined.

An ideal $J \subset  \aa$ is called a {\it reduction ideal} of $\aa$ with respect to $M$ if $J\aa^nM = \aa^{n+1}M$ for some $n\gg 0$. We call  $J$ a {\it minimal reduction} of $\aa$ if  $J$ is a reduction of $\aa$ and $J$ itself does  not contain a proper reduction of $\aa$  (see  \cite{N}). Since $\aa$ is assumed to be an $\mm$-primary ideal, a minimal reduction is  generated by a s.o.p. of $M$. Moreover, we may assume that $J=( x_1, x_2,...,x_d)$ such that the initial forms $x_1^*, ..., x_d^*$ of $x_1,...,x_d$ in $G$ form a filter regular sequence of $G_{\aa}(M)$, that is
$$[(x_1^*,...,x_{i-1}^*)G_{\aa}(M) : x_i^*]_n = [(x_1^*,...,x_{i-1}^*)G_{\aa}(M)]_n,$$
for all $n\gg 0$ and all $i\le d$, see e.g. \cite[Lemma 3.1]{T2}.  

The reduction number of $\aa$ with respect to $J$ and $M$ is the number
 $$r_{J,M}(\aa) = \min\{n\geq 0\mid \aa^{n+1}M=J \aa^n M\}.$$
The reduction number of $\aa$  with respect to  $M$  is the number
$$r_{M}(\aa) = \min\{r_{J,M}(\aa)  \geq 0\mid      J \text{ is a minimal reduction of } \aa  \text{ with respect to }  M  \}.$$
In  the case $M=A$, we just write $r_{J,M}(\aa)$ and $r_{M}(\aa)$ by $r_J(\aa)$ and $r(\aa)$, respectively. We also write $r(A): = r_A(\mm)$  and  call it the reduction number of $A$. Note that $r_{J,M}(\aa) \le r_J(\aa)$ and $r_J(\aa)\le r(\aa)$.

 Let  $L := H_{\mm}^0(M)$ and $\overline{M} = M/L.$ We have 
$$G_\aa(\overline{M}) = \bigoplus_{n\geq 0}(\aa^nM + L)/(\aa^{n+1}M + L)
\cong \bigoplus_{n\geq 0}\aa^nM/(\aa^{n+1}M + \aa^nM\cap L).$$
The relationship between $\reg(G_\aa(M))$ and $\reg(G_\aa(\overline{M})) $  is given by the following
lemma.
\begin{Lemma}\label{a2.2.11 (i)} {\rm (\cite[Lemma 4.3]{L})}
 Let  $L := H_{\mm}^0(M)$ and $\overline{M} = M/L.$ Then,
 $$\reg(G_\aa(M)) \leq \reg(G_\aa(\overline{M})) + \ell(L).$$
\end{Lemma}
In the case $M$ is a $(k,d)$-Buchsbaum with $k>0$, we can replace $ \ell(L)$ by $k$ and some other invariants of $M$ as follows.
\begin{Lemma} \label{a2.2.11 (ii)}
Let $A$ be a Noetherian local ring with an infinite residue field and  $M$  a $(k,d)$-Buchsbaum module with  $k>0$. Then
\begin{itemize}
\item[{\rm i)}] $\reg(G_\aa(M))\leq \reg(G_\aa(\overline{M}))+r_{Q,M}(\aa)+(k-1)d+1,$ where $Q$   is a minimal reduction ideal of $\aa$ with respect to $M.$
\item[{\rm ii)}] $\reg(G_\aa(M))\leq \reg(G_\aa(\overline{M}))+r(A)+(k-1)d+1$. 
\end{itemize}
\end{Lemma}
\begin{pf} i)   There is a short exact sequence 
$$ 0\longrightarrow K\longrightarrow G_\aa(M)\longrightarrow G_\aa(\overline{M})\longrightarrow 0,$$
where 
$$K = \bigoplus_{n\geq 0}(\aa^{n+1}M + \aa^nM\cap L)/\aa^{n+1}M \cong
\bigoplus_{n\geq 0}(\aa^nM \cap  L)/(\aa^{n+1}M\cap L).$$
By  \cite [Corollary 20.19 (b)]{E}, we have 
\begin{equation}\label{eq:1}
\reg(G_\aa(M))\leq \max\{\reg(G_\aa(\overline{M})), \reg(K)\}\leq \reg(G_\aa(\overline{M}))+\reg(K).
\end{equation}
For short, put $r = r_{Q,M}(\aa)$. By Lemma  \ref{a1.3.3}, for all $m\ge 0$, 
  $$\aa^{r+(k-1)d+2 +m}M\cap L\subseteq Q^{(k-1)d+1}M\cap L=0,$$
which yields $K_n = 0$  for all $n\ge r+(k-1)d+2$. Moreover $\ell(K_n) \le \ell(L) < \infty $ for all $n$. Hence
  $$\reg(K) = \max\{n\mid K_n\ne 0\} \leq r+(k-1)d + 1.$$
Using (\ref{eq:1}), we get
$$\reg(G_\aa(M))\leq \reg(G_\aa(\overline{M}))+r+(k-1)d + 1.$$

ii) Let $r' = r(A)$. The proof is similar to that of  i) by noticing the following:
 $$\aa^{r'+(k-1)d+2}M\cap L\subseteq \mm^{r'+(k-1)d+2}M\cap L \subseteq  \qq^{(k-1)d+1}M\cap L=0,$$
 where $\qq$ is a minimal reduction of $\mm$.
\end{pf}

\begin{Lemma}\label{a2.2.12}
Assume that  $M$ is finite $A$-module and $Q$   is a minimal reduction ideal of $\aa$ with respect to $M$. Let $r=r_{Q,M}(\aa).$ Then,
$$\reg(G_\aa(M))\leq \reg(G_Q(\aa^rM))+r.$$
\end{Lemma}
\begin{pf}
There is nothing to prove if  $r=0.$ 
Let  $r\geq 1.$ Note that
$$G_\aa(M)=M/\aa M\oplus \cdots \oplus \aa^{r-1}M/\aa^{r}M\oplus G_Q(\aa^{r}M).$$
Since 
$\ell(\aa^iM/\aa^{i+1}M)<\infty,$ for all $i\leq r-1,$ we get
\begin{equation}\label{eq:11}
 \reg(M/\aa M\oplus \cdots \oplus \aa^{r-1}M/\aa^{r}M)\leq r -1.
\end{equation}
There is a short exact sequence 
$$ 0\longrightarrow G_Q(\aa^{r}M)(-r) \longrightarrow G_\aa(M)\longrightarrow M/\aa M\oplus \cdots \oplus \aa^{r-1}M/\aa^{r}M \longrightarrow 0.$$
By  \cite [Corollary 20.19 (b)]{E}, we have 
\begin{eqnarray}\nonumber
\reg(G_\aa(M)) & \leq & \max\{\reg( G_Q(\aa^{r}M)(-r)),  \reg(M/\aa M\oplus \cdots \oplus \aa^{r-1}M/\aa^{r}M) \} \\
& =  &  \max\{ \reg( G_Q(\aa^{r}M)) +r ,  \reg(M/\aa M\oplus \cdots \oplus \aa^{r-1}M/\aa^{r}M)  \}  \nonumber \\
& \leq  &  \max\{ \reg( G_Q(\aa^{r}M)) +r , r-1 \} \text{ (By (\ref{eq:11}))} \nonumber \\
& = &   \reg(G_Q(\aa^r M)) + r.  \nonumber
\end{eqnarray}
\end{pf} 

\section{Main results}

Recall that $(A,\mm)$ is assumed to be a  Noetherian local ring and $M$ a finitely generated $(k,d)$-Buchsbaum module over $A$, where  $d = \dim A = \dim M$ and  $\aa$ an $\mm$-primary ideal. The purpose of this section is to give bounds on the Castelnuovo-Mumford regularity  $\reg(G_\aa(M))$ in terms of $k$, $d$ and reduction number. 

\begin{Theorem}\label{a3.1.2}
Let $M$ be a  $(k,d)$-Buchsbaum module  with $k\ge 0$ and $Q$ a minimal reduction ideal of $\aa$ with respect to $M$. Set $r:=r_{Q,M}(\aa)$. Then,
\begin{itemize}
\item[{\rm (i) }] $ \reg(G_\aa(M))\leq 2r+k$ {\rm  if } $d=1,$
\item[{\rm (ii) }] $ \reg(G_\aa(M))\leq (2r+k+1)(\ell(M / Q M)+1) + r + 2k-1$ {\rm  if } $d=2,$
\item[{\rm (iii) }] $  \reg(G_\aa(M)) \leq (2r+k+1)^{(d-1)!}[\ell(M / Q M)+1]^{(2d-3)(d-2)!} -1 \mbox{  \rm  if  } d\geq 3.$
\end{itemize}
\end{Theorem}

\begin{pf}  First, assume that $d=1$.  Then  $\overline{M} := M/L$ is  a Cohen-Macaulay module, where $L := H^0_{\mm}(M)$. By Lemma \ref{a2.2.11 (ii)} (i),
$$ \reg(G_\aa (M)) \leq \reg(G_\aa(\overline{M}))+r+k.$$
By Lemma \ref{a2.2.12}, we have
$$\reg(G_\aa(\overline{M}))\leq \reg(G_Q(\aa^r\overline{M}))+r.$$
 Note that   Since $\aa^r\overline{M}\subset  \overline{M}$ is  a Cohen-Macaulay module, $Q$  is generated by a regular element $x$  of $\aa^r\overline{M}$. Hence 
  $G_Q(\aa^rM) \cong \aa^r M/ Q\aa^rM[x]$, which implies $\reg(G_Q(\aa^rM))=0$. Hence $\reg(G_\aa(M))\leq 2r+k$.

Now assume that $d\geq 2.$  We may assume that $Q=(x := x_1, x_2,...,x_d)$ such that the initial forms $x_1^*, ..., x_d^*$ of $x_1,...,x_d$ in $G:= G_{\aa}(A)$ form a filter regular sequence of $G_{\aa}(M)$.  Note that  $\overline{M}$ is also a $(k,d)$-Buchsbaum module, $N:=M/xM$ is  a $(k,d-1)$-Buchsbaum module and $r_{Q',M/xM}(\aa)\leq r,$ $\ell(N/Q'N) \le \ell:= \ell(M/QM)$, where  $Q'=(x_2,...x_d).$

Let   $m\ge \reg(G_\aa(N)) $ be any integer.  Denote by $p_{G_\aa(\overline{M})}$ the Hilbert polynomial of $G_\aa(\overline{M})$, i.e. $p_{G_\aa(\overline{M})}(n) = \ell( G_\aa(\overline{M})_n) = \ell(\aa^n\overline{M}/ \aa^{n+1}\overline{M})$ for $n\gg 0$. Similarly as in the proof of \cite[Theorem 4.4 (ii)]{L},  
we have
$$\reg^1(G_\aa(\overline{M})) \leq m + p_{G_\aa(\overline{M})}(m).$$
 By  \cite[Lemma 3.5]{L} and \cite[Lemma 1.7 (i)]{DH1}, we have
\begin{align*}
p_{G_\aa (\overline{M})}(m) &  \leq  \ell(N/\aa^{m+1}N)   \leq  \ell(N/Q{'}^{m+1}N)  \\
&\leq \ell(N/Q'N)\dbinom{m+d-1}{d-1} \leq \ell \dbinom{m+d-1}{d-1} .
\end{align*}
Using  Lemma \ref{a2.2.11 (ii)} this implies
\begin{eqnarray}\nonumber
\reg(G_\aa(M) & \leq& \reg(G_\aa (\overline{M})) +r + kd \\
& = & \reg^1(G_\aa (\overline{M})) +r + kd \ \   \mbox{(By \cite[Corollary 5.3]{H2})} \nonumber\\
& \leq  & m + \ell \dbinom{m+d-1}{d-1}+r + kd. \label{Induction}
\end{eqnarray}

\noindent If $d= 2,$ then we can conclude that
$$\reg(G_\aa(M)) \leq m+\ell (m+1)+r + k d .$$
By (i) of the theorem one can take  $m = 2r+k.$ Then
\begin{align*}
\reg(G_\aa(M)) &\leq  2r+k +
\ell( 2r+k +1)   + r+2k \\
& = (2r+k+1)(\ell +1) + r +2k-1.
\end{align*}

Finally, let $d\geq 3. $  If $\ell= \ell(M/QM)=1$,  $M$ is a Cohen-Macaulay module and $QM=\mm M=\aa M$. Then one can replace $\aa$ by $Q$ and $Q$ is generated by a regular sequence $Q=(x_1,...,x_d)$. In this case $G_{\aa}(M) \cong M/QM[x_1,...,x_d]$, which yields  $\reg(G_\aa(M)) =0.$ 
So, we can assume that $\ell >1$. 

If $d =3$. By the case $d=2$,  $m \leq (2r+k+1)(\ell +1) + r +2k-1$.
Hence
\begin{align*}
 \reg(G_\aa(M)) & \leq m + \frac{(m+1)(m+2)}{2} +r+3k  \\ 
& \leq (2r+k+1)^2 (\ell+1)^3-1.
\end{align*}

If $d \ge 4$, by the inductive assumption we can take
$$ m = (2r+k+1)^{(d-2)!}(\ell +1)^{(2d-5)(d-3)!} -1 .$$
Since $\ell > 1$, we have  $m\geq 2$. This implies
\begin{align*}
\reg(G_\aa(M)) & \leq m +  \ell(m+1)^{d-1}+r+kd\\
& = (2r+k+1)^{(d-1)!} \ell \left(\ell+1 \right)^{(2d-5)(d-3)! (d-1)}+r+dk \\
&\leq (2r+k+1)^{(d-1)!}  \left(\ell+1 \right)^{(2d-5)(d-3)! (d-1)+1}-1.
\end{align*}
Since 
$  (2d-5)(d-3)! (d-1)+1 \le (2d-3)(d-2)!$
for all $d \ge 3$, the above inequality yields
$$ \reg(G_\aa(M))  \le   (2r+k+1)^{(d-1)!} \left(\ell+1 \right)^{ (2d-3)(d-2)!}-1,$$
as required.
\end{pf}

In the next result, we replace $r_{Q,M}(\aa)$ by the reduction number  $r(A)$ of $A$, but we need to add  the multiplicity $e(\aa,M)$ of $M$ with respect to $\aa$.

\begin{Theorem}\label{a3.1.22}
Let $M$ be a $(k,d)$-Buchsbaum module  with $k\ge 0$, $Q$ a minimal reduction ideal of $\aa$ with respect to $M$ and $r(A)$ reduction number of $A$. Then,
\begin{itemize}
\item[{\rm (i) }] $ \reg(G_\aa(M))\leq e(\aa,M)+ r(A)+k -1$ {\rm  if } $d=1,$
\item[{\rm (ii) }] $ \reg(G_\aa(M))\leq (e(\aa,M)+ r(A)+k)(\ell(M / Q M)+1) + r(A) + 2k-1$ {\rm  if } $d=2,$
\item[{\rm (iii) }] $  \reg(G_\aa(M)) \leq [e(\aa,M)+ r(A)+ k]^{(d-1)!}[\ell(M / Q M)+1]^{(2d-3)(d-2)!} -1 \mbox{  \rm  if  } d\geq 3.$
\end{itemize}
\end{Theorem}
\begin{pf}
If $d=1$ then  $\overline{M}$ is  a Cohen-Macaulay module. By Lemma \ref{a2.2.11 (ii)} ii),
$$ \reg(G_\aa (M)) \leq \reg(G_\aa(\overline{M}))+r(A)+k.$$
By \cite[Corollary 5.3]{H2} and \cite[Lemma 3.1]{HM}, we have 
$$\reg(G_\aa(\overline{M})) = \reg^1(G_\aa(\overline{M})) \le e(\aa,\overline{M}) -1 = e(\aa,M) -1.$$
Hence 
$$\reg(G_\aa(M))\leq e(\aa,M)+ r(A)+k -1.$$

If $d\geq 2$,   the proof is similar to that of Theorem \ref{a3.1.2}, so that we omit its detail and give here only the sketch.  In this case, we also have the induction formula (\ref{Induction}).  For $d\geq 3 $,  we can also restrict to the case $\ell >1$. Then one should first check the case $d=3$, and do by induction for $d\ge 4$.
 \end{pf}

If $M$ is a Cohen-Macaulay module, $\reg(G_{\aa}(M))$ can be bounded in terms of  the multiplicity $e(\aa,M)$, see  \cite[Corollary 3.4]{RTV} and \cite[Corollary 4.5]{L}. The following result gives a different bound which  involves the reduction number. It is easy to give examples when $r_M(\aa)$ is rather small compared with $e(\aa,M)$. So the new bound some times give a slightly better bound than that of  \cite[Corollary 3.4]{RTV} and \cite[Corollary 4.5]{L}.

\begin{Corollary}  \label{aa}  Let $M$ be a Cohen-Macaulay module. Then
  $$\reg(G_{\aa}(M))\leq \begin{cases} 
2r_M(\aa) &  \text{ if } d=1,\\
(2r_M(\aa)+1)(e(\aa,M)+1) + r_M(\aa)  -1 &  \text{ if } d=2,\\
(2r_M(\aa)+1)^{(d-1)!}[e(\aa,M)+1]^{(2d-3)(d-2)!} -1 &  \text{ if } d \ge 3.
\end{cases}$$
\end{Corollary}
\begin{pf}
 The Cohen-Macaulay property means $M$ is a $(0,d)$-Buchsbaum module, i.e. one can take $k=0$.  In this case, we need to take a minimal reduction $Q$ such that $r_{Q,M}(\aa) = r_M(\aa)$. Then $\ell(M/QM) = e(Q,M) = e(\aa,M)$, and the claim immediately follows from  Theorem \ref{a3.1.2}.
\end{pf}

\begin{Corollary} Let $M$ be a Buchsbaum module and $Q$ a minimal reduction of $\aa$ with respect to $M$. Then
\begin{itemize}
\item[{ \rm (i)}]  $\reg(G_{\aa}(M))\leq \begin{cases} 
2r_{Q,M}(\aa)+1 &  \text{ if } d=1,\\
(2r_{Q,M}(\aa)+2)(\ell(M / Q M)+1) + r_{Q,M}(\aa) + 1 &  \text{ if } d=2,\\
(2r_{Q,M}(\aa)+2)^{(d-1)!}[\ell(M / Q M)+1]^{(2d-3)(d-2)!} -1 &  \text{ if } d \ge 3.
\end{cases}$

\item[{\rm (ii)}]  $\reg(G_{\aa} (M))\leq \begin{cases} 
e(\aa,M) + r(A) &  \text{ if } d=1,\\
 (e(\aa,M)+ r(A) +1)[\ell(M/QM)+1] + r(A) + 1&  \text{ if } d=2,\\
(e(\aa,M)+ r(A) + 1)^{(d-1)!}[\ell(M/QM)+1]^{(2d-3)(d-2)!} -1&  \text{ if } d \ge 3.
\end{cases}$
\end{itemize}
\end{Corollary}
\begin{pf}
 A Buchsbaum module is a $(1,d)$-Buchsbaum module. Therefore, applying  Theorem  \ref{a3.1.2} and  Theorem  \ref{a3.1.22} to $k=1$, we  get the claims.
\end{pf}

Using also Lemma \ref{a1.3.3}, we can improve a little bit Theorem \ref{a3.1.2} and Theorem \ref{a3.1.22}. Let $k' = \lceil \frac{k}{o(\aa)} \rceil$, where  $o(\aa) = \max\{t | \aa \subseteq \mm^t\}$ is the order of the ideal $\aa$.

\begin{Corollary}\label{c1}
Let $M$ be a  $(k,d)$-Buchsbaum module  with $k\ge 0$ and $Q $   a minimal reduction ideal of $\aa$ with respect to $M$. Then
\begin{itemize}
\item[{\rm (i)}]  $\reg(G_{\aa}(M))\leq \begin{cases} 
2r+k' &  \text{ if } d=1,\\
(2r+k'+1)(\ell(M / Q M)+1) + r + 2k'-1 &  \text{ if } d=2,\\
(2r+k'+1)^{(d-1)!}[\ell(M / Q M)+1]^{(2d-3)(d-2)!} -1 &  \text{ if } d \ge 3,
\end{cases}$
where $r=r_{Q,M}(\aa) $.
\item[{\rm (ii)}]  $\reg(G_{\aa}(M))\leq \begin{cases} 
e+ r'+k' -1&  \text{ if } d=1,\\
 (e+ r'+k')[\ell(M/QM)+1] + r' + 2k'-1&  \text{ if } d=2,\\
(e+ r'+ k')^{(d-1)!}[\ell(M/QM)+1]^{(2d-3)(d-2)!} -1&  \text{ if } d \ge 3,
\end{cases}$
where $e: = e(\aa,M)$ and $r' = r(A)$.
\end{itemize}
\end{Corollary}

\begin{pf} Assume that $Q= (x_1,...,x_d)$. Then $x_1,...,x_d \in \aa \subseteq m^{o(\aa)}$. By Lemma   \ref{a1.3.3}, we have
$$Q^{(k'-1)d+1}M\cap H_{\mm}^0(M)=0.$$
Now, replacing $k$ by $k'$ in the proofs of Theorem \ref{a3.1.2} and  Theorem \ref{a3.1.22},  we get the claims.
\end{pf}
In particular, we get

\begin{Corollary} Let $M$ be n $(k,d)$-Buchsbaum module with $k > 0$. Assume that $\aa \subseteq \mm^k$. Let $Q $ be  a minimal reduction ideal of $\aa$ with respect to $M$. Then
\begin{itemize}
\item[{ \rm (i)}]  $\reg(G_{\aa}(M))\leq \begin{cases} 
2r+1 &  \text{ if } d=1,\\
(2r+2)(\ell(M / Q M)+1) + r + 1 &  \text{ if } d=2,\\
(2r+2)^{(d-1)!}[\ell(M / Q M)+1]^{(2d-3)(d-2)!} -1 &  \text{ if } d \ge 3,
\end{cases}$
where $r=r_{Q,M}(\aa) $.
\item[{\rm (ii)}]  $\reg(G_{\aa}(M))\leq \begin{cases} 
e+ r'&  \text{ if } d=1,\\
 (e+ r'+1)[\ell(M/QM)+1] + r' + 1&  \text{ if } d=2,\\
(e+ r'+ 1)^{(d-1)!}[\ell(M/QM)+1]^{(2d-3)(d-2)!} -1&  \text{ if } d \ge 3,
\end{cases}$
where $e: = e(\aa,M)$ and $r' = r(A)$.
\end{itemize}
\end{Corollary}

Hoa and Vogel in \cite[Corollary 2.3]{HV} showed that a $k$-Buchsbaum modules $M$ is a $(kd,d)$-Buchsbaum modules. From Theorem \ref{a3.1.2} and 
Theorem \ref{a3.1.22}, we then obtain  
\begin{Theorem} \label{main}
Let $M$ be an  $k$-Buchsbaum module  with $k\ge 0$ and $Q $   a minimal reduction ideal of $\aa$ with respect to $M$. Then
\begin{itemize}
\item[{ \rm (i)}]  $\reg(G_{\aa}(M))\leq \begin{cases} 
2r+k &  \text{ if } d=1,\\
(2r+2k+1)(\ell(M / Q M)+1) + r + 4k-1 &  \text{ if } d=2,\\
(2r+kd+1)^{(d-1)!}[\ell(M / Q M)+1]^{(2d-3)(d-2)!} -1 &  \text{ if } d \ge 3,
\end{cases}$
where $r=r_{Q,M}(\aa) $.

\item[{\rm (ii)}]  $\reg(G_{\aa}(M))\leq \begin{cases} 
e+ r'+k -1&  \text{ if } d=1,\\
 (e+ r'+2k)[\ell(M/QM)+1] + r' + 4k-1&  \text{ if } d=2,\\
(e+ r'+ kd)^{(d-1)!}[\ell(M/QM)+1]^{(2d-3)(d-2)!} -1&  \text{ if } d \ge 3,
\end{cases}$
where $e: = e(\aa,M)$ and $r' = r(A)$.
\end{itemize}
\end{Theorem}

 \vskip0.5cm

\noindent {\bf Acknowledgment}:  The author is grateful to Prof. L. T. Hoa for  many discussions and comments. The author is partially supported by Project B2017-HDT-05.


\begin{thebibliography}{10}

\bibitem{CST} 
N. T. Cuong, P. Schenzel, N. V. Trung, Verallgemeinerte Cohen-Macaulay-Moduln, {\it Math. Nachr.} {\bf 85} (1978) 57-73 .

\bibitem{DH1} L. X. Dung and L. T. Hoa,    Castelnuovo-Mumford Regularity of Associated Graded Modules and Fiber Cones of Filtered Modules, {\it Comm. Algebra} {\bf 40}  (2012) 404-422.

\bibitem{DH2} L. X. Dung and L. T. Hoa,   Dependence of Hilbert coefficients,  {\it manuscripta math.} {\bf 149} (2016) 404-422.

\bibitem{DH3} L. X. Dung and L. T. Hoa,   A note on Castelnuovo-Mumford regularity  and   Hilbert coefficients,  {\it J. Algebra Appl.} (2019), to appear, doi: 10.1142/S0219498819501913.

 \bibitem{E} D. Eisenbud, {\it Commutative Algebra with a View toward Algebraic Geometry}, Springer-Verlag, 1995.


 \bibitem{H2}
 L.T. Hoa, Reduction numbers of equimultiple ideal, {\it J. Pure Appl. Algebra} {\bf 109} (1996) 111-126.
 
\bibitem {HM}
L.T. Hoa and C. Miyazaki,   Bounds on Castelnuovo-Mumford
regularity for generalized Cohen-Macaulay graded rings, {\it Math.
Ann.} {\bf 301} (1995) 587-598.

\bibitem {HV}
L.T. Hoa and W. Vogel,  Castelnuovo-Mumford regularity and hyperplane sections,
{\it J. Algebra} {\bf 163} (1994) 348-365.

\bibitem {L}
C.H. Linh, Upper bound for Castelnuovo-Mumford
regularity of associated graded modules, {\it Comm. Algebra} {\bf 33}(2005) 1817-1831.

\bibitem{N}
D.G. Northcott and D. Rees, Reductions of ideals in local rings, {\it Proc. Camb. Phil. Soc.} {\bf 50} (1954) 145-158.
 \bibitem {RTV}
M.E. Rossi, N.V. Trung and G. Valla,  Castelnuovo-Mumford
regularity and extended degree, {\it Trans. Amer. Math. Soc.} {\bf
355} (2003)  1773-1786.

\bibitem {T1}
N.V. Trung,  Towards a theory of generalized Cohen-Macalay
modules, {\it Nagoya Math. J.}  {\bf 102} (1986) 1-49.

\bibitem {T2}
N.V. Trung,  Reduction exponent and degree bound for the
defining equations of graded rings, {\it Proc. Amer. Math. Soc.} {\bf 101} (1987) 229-236.

\bibitem {Va}
W. V. Vasconcelos, Computational methods in commutative algebra and algebraic geometry. With chapters by D. Eisenbud, D.R. Grayson, J. Herzog and M. Stillman. Algorithms and Computation in Math. 2 (Springer-Verlag, Berlin, 1998).


\end{thebibliography}
\end{document}